# The Score, Accuracy, and Certainty Functions determine a Total Order on the Set of Neutrosophic Triplets (T, I, F)

Florentin Smarandache
Division of Mathematics, Physical and Natural Sciences
University of New Mexico, 705 Gurley Ave., Gallup, NM 87301, USA

**Abstract**

In this paper we prove that the Single-Valued (and respectively Interval-Valued, as well as Subset-Valued) Score, Accuracy, and Certainty Functions determine a total order on the set of neutrosophic triplets (T, I, F). This total order is needed in the neutrosophic decision-making applications.

## 1. Introduction

We reveal the easiest to use single-valued neutrosophic score, accuracy, and certainty functions that exist in the literature and the algorithm how to use them all together. We present Xu and Da's Possibility Degree that an interval is greater than or equal to another interval, and we prove that this method is equivalent to the intervals' midpoints comparison. Also, Hong-yu Zhang et al.'s interval-valued neutrosophic score, accuracy, and certainty functions are listed, that we simplify these functions. Numerical examples are provided.

## 2. Single-Valued Neutrosophic Score, Accuracy, and Certainty Functions

We firstly present the most known and used in literature single-valued score, accuracy, and certainty functions.

Let M be the set of single-valued neutrosophic triplet numbers,

$M = \{(T, I, F), \text{where } T, I, F \in [0, 1], 0 \leq T + I + F \leq 3\}$.

Let $N = (T, I, F) \in M$ be a generic single-valued neutrosophic triplet number. Then:

$T$ = truth (or membership) represents the positive quality of $N$;

$I$ = indeterminacy represents a negative quality of $N$,

hence $1 - I$ represents a positive quality of $N$;

$F$ = falsehood (or nonmembership) represents also a negative quality of $N$, hence $1 - F$ represents a positive quality of $N$.

We present the three most used and best functions in the literature:

### 2.1) The Single-Valued Neutrosophic Score Function

$$s: M \to [0, 1]$$

$$s(T, I, F) = \frac{T + (1 - I) + (1 - F)}{3} = \frac{2 + T - I - F}{3}$$

that represents the average of positiveness of the single-valued neutrosophic components $T$, $I$, and $F$.

### 2.2) The Single-Valued Neutrosophic Accuracy Function

$$a: M \to [-1, 1]$$

$$a(T, I, F) = T - F$$

### 2.3) The Single-Valued Neutrosophic Certainty Function

$$c: M \to [0, 1]$$

$$c(T, I, F) = T$$

## 3. Algorithm for Ranking the Single-Valued Neutrosophic Triplets

Let $(T_1, I_1, F_1)$ and $(T_2, I_2, F_2)$ be two single-valued neutrosophic triplets from $M$, i.e. $T_1, I_1, F_1, T_2, I_2, F_2 \in [0, 1]$.

Apply the *Neutrosophic Score Function*.

1. If $s(T_1, I_1, F_1) > s(T_2, I_2, F_2)$, then $(T_1, I_1, F_1) > (T_2, I_2, F_2)$.

2. If $s(T_1, I_1, F_1) < s(T_2, I_2, F_2)$, then $(T_1, I_1, F_1) < (T_2, I_2, F_2)$.

3. If $s(T_1, I_1, F_1) = s(T_2, I_2, F_2)$, then apply the *Neutrosophic Accuracy Function*:

   3.1 If $a(T_1, I_1, F_1) > a(T_2, I_2, F_2)$, then $(T_1, I_1, F_1) > (T_2, I_2, F_2)$.

   3.2 If $a(T_1, I_1, F_1) < a(T_2, I_2, F_2)$, then $(T_1, I_1, F_1) < (T_2, I_2, F_2)$.

   3.3 If $a(T_1, I_1, F_1) = a(T_2, I_2, F_2)$, then apply the *Neutrosophic Certainty Function*.

   3.3.1 If $c(T_1, I_1, F_1) > c(T_2, I_2, F_2)$, then $(T_1, I_1, F_1) > (T_2, I_2, F_2)$.

   3.3.2 If $c(T_1, I_1, F_1) < c(T_2, I_2, F_2)$, then $(T_1, I_1, F_1) < (T_2, I_2, F_2)$.

3.3.1 If $c(T_1, I_1, F_1) = c(T_2, I_2, F_2)$, then $(T_1, I_1, F_1) \equiv (T_2, I_2, F_2)$, i.e. $T_1 = T_2$, $I_1 = I_2$, $F_1 = F_2$.

**3.1. Theorem**

We prove that the single-valued neutrosophic score, accuracy, and certainty functions all together form a *total order relationship* on $M$. Or:

for any two single-valued neutrosophic triplets $(T_1, I_1, F_1)$ and $(T_2, I_2, F_2)$ we have:

a) Either $(T_1, I_1, F_1) > (T_2, I_2, F_2)$

b) Or $(T_1, I_1, F_1) < (T_2, I_2, F_2)$

c) Or $(T_1, I_1, F_1) \equiv (T_2, I_2, F_2)$, which means that $T_1 = T_2$, $I_1 = I_2$, $F_1 = F_2$.

Therefore, on the set of single-valued neutrsophic triplets $M = \{(T, I, F), \text{with } T, I, F \in [0, 1], 0 \leq T + I + F \leq 3\}$, the score, accuracy, and certainty functions altogether form a *total order relationship*.

***Proof.***

Firstly we apply the score function.

The only problematic case is when we get equality:

$$s(T_1, I_1, F_1) = s(T_2, I_2, F_2)$$

That means:

$$\frac{2 + T_1 - I_1 - F_1}{3} = \frac{2 + T_2 - I_2 - F_2}{3}$$

or $T_1 - I_1 - F_1 = T_2 - I_2 - F_2$.

Secondly we apply the accuracy function.

Again the only problematic case is when we get equality:

$$a(T_1, I_1, F_1) = a(T_2, I_2, F_2)$$

or $T_1 - F_1 = T_2 - F_2$.

Thirdly, we apply the certainty function.

Similarly, the only problematic case may be when we get equality:

$$c(T_1, I_1, F_1) = c(T_2, I_2, F_2)$$

or $T_1 = T_2$.

For the most problematic case, we got the following linear algebraic system of 3 equations of 6 variables:

$$\begin{cases} T_1 - I_1 - F_1 = T_2 - I_2 - F_2 \\ T_1 - F_1 = T_2 - F_2 \\ T_1 = T_2 \end{cases}$$

Let's solve it.

Since $T_1 = T_2$, replacing this into the second equation we get $F_1 = F_2$.

Now, replacing both $T_1 = T_2$ and $F_1 = F_2$ into the first equation, we get $I_1 = I_2$.

Therefore the two neutrosophic triplets are identical: $(T_1, I_1, F_1) \equiv (T_2, I_2, F_2)$, i.e. equivalent (or equal), or $T_1 = T_2$, $I_1 = I_2$, and $F_1 = F_2$.

In conclusion, for any two single-valued neutrosophic triplets, either one is bigger than the other, or both are equal (identical).

## 4. Definition of Neutrosophic Negative Score Function

We have introduce in 2017 for the first time [1] the *Average Negative Quality Neutrosophic Function* of a single-valued neutrosophic triplet, defined as: (2)

$$s^- : [0,1]^3 \to [0,1], \; s^-(t,i,f) = \frac{(1-t)+i+f}{3} = \frac{1-t+i+f}{3}.$$

### 4.1. Theorem

The average positive quality (score) neutrosophic function and the average negative quality neutrosophic function are complementary to each other, or

$$s^+(t,i,f) + s^-(t,i,f) = 1. \qquad (3)$$

*Proof.*

$$s^+(t,i,f) + s^-(t,i,f) = \frac{2+t-i-f}{3} + \frac{1-t+i+f}{3} = 1. \qquad (4)$$

The *Neutrosophic Accuracy Function* has been defined by:

h: [0, 1]³ → [-1, 1], h(t, i, f) = t - f.  (5)

We have also introduce [1] for the first time the *Extended Accuracy Neutrosophic Function*, defined as follows:

$h_e$: [0, 1]³ → [-2, 1], $h_e$(t, i, f) = t − i − f,  (6)

which varies on a range: from the worst negative quality (-2) [or minimum value], to the best positive quality (+1) [or maximum value].

### 4.2. Theorem

If $s(T_1, I_1, F_1) = s(T_2, I_2, F_2)$, $a(T_1, I_1, F_1) = a(T_2, I_2, F_2)$, and $c(T_1, I_1, F_1) = c(T_2, I_2, F_2)$, then $T_1 = T_2$, $I_1 = I_2$, $F_1 = F_2$, or the two neutrosophic triplets are identical: $(T_1, I_1, F_1) \equiv (T_2, I_2, F_2)$.

*Proof:*
It results from the proof of Theorem 3.1.

## 5. Xu and Da's Possibility Degree

Xu and Da [3] have defined in 2002 the possibility degree $P(.)$ that an interval is greater than another interval:

$$[a_1, a_2] \geq [b_1, b_2]$$

for $a_1, a_2, b_1, b_2 \in [0,1]$ and $a_1 \leq a_2, b_1 \leq b_2$, in the following way:

$$P([a_1, a_2] \geq [b_1, b_2]) = max\left\{1 - max\left(\frac{b_2 - a_1}{a_2 - a_1 + b_2 - b_1}, 0\right), 0\right\},$$

where $a_2 - a_1 + b_2 - b_1 \neq 0$ (i.e. $a_2 \neq a_1$ or $b_2 \neq b_1$.

They proved the following:

### 5.1. Properties

1) $P([a_1, a_2] \geq [b_1, b_2]) \in [0,1]$;
2) $P([a_1, a_2] \approx [b_1, b_2]) = 0.5$;
3) $P([a_1, a_2] \geq [b_1, b_2]) + P([b_1, b_2] \geq [a_1, a_2]) = 1$.

### 5.2. Example

Let $[0.4, 0.7]$ and $[0.3, 0.6]$ be two intervals.
Then,

$$P([0.4, 0.7] \geq [0.3, 0.6]) = max\left\{1 - max\left(\frac{0.6 - 0.4}{0.7 - 0.4 + 0.6 - 0.3}, 0\right), 0\right\}$$

$$= max\left\{1 - max\left(\frac{0.2}{0.6}\right), 0\right\} = max\left\{1 - \frac{0.2}{0.6}, 0\right\} = \frac{0.4}{0.6} \approx 0.66 > 0.50,$$

therefore $[0.4, 0.7] \geq [0.3, 0.6]$.

The opposite:

$$P\big((0.3, 0.6) \geq ([0.4, 0.7])\big) = max\left\{1 - max\left(\frac{0.7 - 0.3}{0.6 - 0.3 + 0.7 - 0.4}, 0\right), 0\right\}$$

$$= max\left\{1 - max\left(\frac{0.4}{0.6}, 0\right), 0\right\} = max\left\{1 - \frac{0.4}{0.6}, 0\right\} = \frac{0.2}{0.6} \approx 0.33 < 0.50,$$

therefore $[0.3, 0.6] \leq [0.4, 0.7]$.

We see that

$$P([0.4, 0.7] \geq [0.3, 0.6]) + P([0.3, 0.6] \geq [0.4, 0.7]) = \frac{0.4}{0.6} + \frac{0.2}{0.6} = 1.$$

Another method of ranking two intervals is the midpoint one.

## 6. Midpoint Method

Let $A = [a_1, a_2]$ and $B = [b_1, b_2]$ be two intervals included in or equal to [0, 1], with $m_A = (a_1 + a_2)/2$ and $m_B = (b_1 + b_2)/2$ the midpoints of $A$ and respectively $B$. Then:

1) If $m_A < m_B$ then $A < B$.

2) If $m_A > m_B$ then $A > B$.

3) If $m_A = m_B$ then $A =_N B$, i.e. A is neutrosophically equal to B.

### 6.1 Example

1) We take the previous example,

where $A = [0.4, 0.7]$, and $m_A = \frac{0.4 + 0.7}{2} = 0.55$;

and $B = [0.3, 0.6]$, and $m_B = \frac{0.3 + 0.6}{2} = 0.45$.

Since $m_A = 0.55 > 0.45 = m_B$, we have $A > B$.

Let $C = [0.1, 0.7]$ and $D = [0.3, 0.5]$.

Then $m_C = \frac{0.1 + 0.7}{2} = 0.4$, and $m_D = \frac{0.3 + 0.5}{2} = 0.4$.

Since $m_C = m_D = 0.4$, we get $C =_N D$.

Let's verify the ranking relationship between $C$ and $D$ using Xu and Da's possibility degree method.

$$P([0.1, 0.7] \geq [0.3, 0.5]) = max\left\{1 - max\left(\frac{0.5 - 0.1}{0.7 - 0.1 + 0.5 - 0.3}, 0\right), 0\right\}$$

$$= max\left\{1 - max\left(\frac{0.4}{0.8}, 0\right), 0\right\} = max\left\{1 - \frac{0.4}{0.8}, 0\right\} = max\left\{\frac{0.4}{0.8}, 0\right\} = 0.5;$$

and $P([0.3, 0.5] \geq [0.1, 0.7]) = max\left\{1 - max\left(\frac{0.7-0.3}{0.5-0.3+0.7-0.1}, 0\right), 0\right\} = max\{1 - max\left(\frac{0.4}{0.8}, 0\right), 0\} = max\{1 - \frac{0.4}{0.8}, 0\} = max\{\frac{0.4}{0.8}, 0\} = 0.5;$

thus, $[0.1, 0.7] =_N [0.3, 0.5]$.

### 6.2. Corollary

The possibility method for two intervals having the same midpoint gives always 0.5.

For example:
*p([0.3, 0.5] ≥ [0.2, 0.6])* = *max{1 - max( ((0.6-0.3) / (0.5-0.3 + 0.6-0.2)), 0 ), 0}* =
= *max{1 - max( ((0.3) / (0.6)), 0 ), 0}*
= *max{1 - max( 0.5, 0 ), 0}* = *0.5*.

Similarly,
*p([0.2, 0.6] ≥ [0.3, 0.5])* = *max{1 - max( ((0.5-0.2) / (0.6-0.2 + 0.5-0.3)), 0 ), 0}* = *0.5*.
Hence, none of the intervals [0.3, 0.5] and [0.2, 0.6]) is bigger than the other.

Therefore, we may consider that the intervals $[0.3, 0.5] =_N [0.2, 0.6]$ are neutrosophically equal (or neutrosophically equivalent).

### 7. Normalized Hamming Distance between Two Intervals

Let's consider the *Normalized Hamming Distance* between two intervals *[a₁, a₂]* and *[b₁, b₂]*
$$h : int([0, 1]) \times int([0, 1]) \rightarrow [0, 1]$$
defined as follows:
$$h([a_1, b_1], [a_2, b_2]) = \tfrac{1}{2}(|a_1 - b_1| + |a_2 - b_2|).$$

### 7.1. Theorem

7.1.1. The Normalized Hamming Distance between two intervals having the same midpoint and the negative-ideal interval *[0, 0]* is the same.
7.1.2. The Normalized Hamming Distance between two intervals having the same midpoint and the positive-ideal interval *[1, 1]* is also the same (Jun Ye [4, 5]).

*Proof.*

Let $A = [m - a, m + a]$ and $B = [m - b, m + b]$ be two intervals from *[0, 1]*, where *m-a, m+a, m-b, m+b, a, b, m* ∈ *[0, 1]*. A and B have the same midpoint *m*.

7.1.1. $h([m - a, m + a], [0, 0]) = ½(|m – a - 0| + |m + a - 0|) = ½(m – a + m + a) = m$, and
$h([m - b, m + b], [0, 0]) = ½(|m – b - 0| + |m + b - 0|) = ½(m – b + m + b) = m$,

7.1.2. $h([m - a, m + a], [1, 1]) = ½(|m – a - 1| + |m + a - 1|) = ½(1 - m + a + 1 - m - a) = 1 - m$, and
$h([m - b, m + b], [1, 1]) = ½(|m – b - 1| + |m + b - 1|) = ½(1 - m + b + 1 - m - b) = 1 – m$.

## 8. Xu and Da's Possibility Degree Method is equivalent to the Midpoint Method

We prove the following:

### 8.1. Theorem

The Xu and Da's Possibility Degree Method is equivalent to the Midpoint Method in ranking two intervals included in *[0, 1]*.

*Proof*

Let *A* and *B* be two intervals included in [0,1]. Without loss of generality, we write each interval in terms of each midpoint:

$A = [m_1 - a, m_1 + a]$ and $B = [m_2 - b, m_2 + b]$,

where $m_1, m_2 \in [0,1]$ are the midpoints of A and respectively B, and $a, b \in [0,1]$, $A, B \subseteq [0,1]$.

(For example, if $A = [0.4, 0.7]$, $m_A = \frac{0.4+0.7}{2} = 0.55$, 0.55-0.4=0.15, then $A = [0.55 - 0.15, 0.55 + 0.15]$).

1) First case: $m_1 < m_2$. According to the Midpoint Method, we get $A < B$. Let's prove the same inequality results with the second method.

Let's apply Xu and Da's Possibility Degree Method:

$P(A \geq B) = P([m_1 - a, m_1 + a] \geq [m_2 - b, m_2 + b])$

$$= max\left\{1 - max\left(\frac{(m_2 + b) - (m_1 - a)}{(m_1 + a) - (m_1 - a) + (m_2 + b) - (m_2 - b)}, 0\right), 0\right\}$$

$$= max\left\{1 - max\left(\frac{m_2 - m_1 + a + b}{2a + 2b}, 0\right), 0\right\}$$

$$= max\left\{1 - \frac{m_2 - m_1 + a + b}{2a + 2b}, 0\right\}, \text{because } m_1 < m_2,$$

$$= max\left\{\frac{2a + 2b - m_2 + m_1 - a - b}{2a + 2b}, 0\right\} = max\left\{\frac{a + b + m_1 - m_2}{2a + 2b}, 0\right\}$$

*i)* If $a + b + m_1 - m_2 \leq 0$, then $p(A \geq B) = max\left\{\frac{a+b+m_1-m_2}{2a+2b}, 0\right\} = 0$, hence $A < B$.

*ii)* If $a + b + m_1 - m_2 > 0$,

then $p(A \geq B) = max\left\{\frac{a+b+m_1-m_2}{2a+2b}, 0\right\} = \frac{a+b+m_1-m_2}{2a+2b} > 0$.

We need to prove that $\frac{a+b+m_1-m_2}{2a+2b} < 0.5$,

or $a + b + m_1 - m_2 < 0.5(2a + 2b)$,

or $a + b + m_1 - m_2 < a + b$,

or $m_1 - m_2 < 0$,

or $m_1 < m_2$, which is true according to the first case assumption.

2) Second case: $m_1 = m_2$. According to the Midpoint Method, $A$ is neutrosophically equal to $B$ (we write $A =_N B$).

Let's prove that we get the same result with Xu and Da's Method.

Then $A = [m_1 - a, m_1 + a]$, and $B = [m_1 - b, m_1 + b]$.

Let's apply Xu and Da's Method:

$$P(A \geq B) = max\left\{1 - max\left(\frac{(m_1 + b) - (m_1 - a)}{(m_1 + a) - (m_1 - a) + (m_1 + b) - (m_1 - b)}, 0\right), 0\right\}$$

$$= max\left\{1 - max\left(\frac{a+b}{2a+2b}, 0\right), 0\right\} = max\left\{1 - \frac{1}{2}, 0\right\} = 0.5$$

Similarly:

$$P(B \geq A) = max\left\{1 - max\left(\frac{(m_1 + a) - (m_1 - b)}{(m_1 + b) - (m_1 - b) + (m_1 + a) - (m_1 - a)}, 0\right), 0\right\}$$

$$= max\left\{1 - max\left(\frac{a+b}{2a+2b}, 0\right), 0\right\} = 0.5$$

Therefore, again $A =_N B$.

3) If $m_1 > m_2$, according to the Midpoint Method, we get $A > B$.

Let's prove the same inequality using Xu and Da's Method.

$$P(A \geq B) = P([m_1 - a, m_1 + a] \geq [m_2 - b, m_2 - b]) =$$

$$max\left\{1 - max\left(\frac{(m_2 + b) - (m_1 - a)}{(m_1 + a) - (m_1 - a) + (m_2 + b) - (m_2 - b)}, 0\right), 0\right\}$$

$$= max\left\{1 - max\left(\frac{m_2 - m_1 + a + b}{2a + 2b}, 0\right), 0\right\}$$

*i)* If $m_2 - m_1 + a + b \leq 0$, then $P(A \geq B) = max\{1 - 0, 0\} = 1$, therefore $A > B$.

*ii)* If $m_2 - m_1 + a + b > 0$, then

$$P(A \geq B) = max\left\{1 - \frac{m_2 - m_1 + a + b}{2a + 2b}, 0\right\} = \frac{2a + 2b - m_2 + m_1 - a - b}{2a + 2b}$$
$$= \frac{a + b + m_1 - m_2}{2a + 2b}$$

We need to prove that $\frac{a+b+m_1-m_2}{2a+2b} > 0.5$,

or $a + b + m_1 - m_2 > 0.5(2a + 2b)$

or $a + b + m_1 - m_2 > a + b$

or $m_1 - m_2 > 0$

or $m_1 > m_2$, which is true according to the third case. Thus $A > B$.

### 8.2. Consequence

All intervals, included in [0, 1], that have the same midpoint are considered neutrosophically equal.

$C(m) = \{[m - a, m + a], \text{where all } m, a, m - a, m + a \in [0, 1]\}$

represents the class of all neutrosophically equal intervals included in [0, 1] whose midpoint is $m$.

*i)* If $m = 0$ or $m = 1$, there is only one interval centered in 0, i.e. [0, 0], and only one interval centered in 1, i.e. [1, 1].

*ii)* If $m \notin \{0, 1\}$, there are infinitely many intervals from [0, 1], centered in $m$.

### 8.3. Consequence

Remarkably we can rank an interval $[a, b] \subseteq [0,1]$ with respect to a number $n \in [0, 1]$ since the number may be transformed into an interval $[n, n]$ as well.

For example $[0.2, 0.8] > 0.4$ since the midpoint of $[0.2, 0.8]$ is 0.5, and the midpoint of $[0.4, 0.4] = 0.4$, hence $0.5 > 0.4$.

Similarly, $0.7 > (0.5, 0.8)$.

## 9. Interval (-Valued) Neutrosophic Score, Accuracy, and Certainty Functions

Let $T, I, F \subseteq [0,1]$ be three open, semi-open / semi-closed, or closed intervals.

Let $T^L = infT$ and $T^U = supT$; $I^L = infI$ and $I^U = supI$; $F^L = infF$ and $F^U = supF$.

Let $T^L, T^U, I^L, I^U, F^L, F^U \in [0,1]$, with $T^L \leq T^U, I^L \leq I^U, F^L \leq F^U$.

We consider all possible types of intervals: open (*a, b*), semi-open / semi-closed (*a, b*] and [*a, b*), and closed [*a, b*]. For simplicity of notations, we are using only [*a, b*], but we understand all types.

Then $A = ([T^L, T^U], [I^L, I^U], [F^L, F^U])$ is an Interval Neutrosophic Triplet.

$T^L$ is the lower limit of the interval $T$,

$T^U$ is the upper limit of the interval $T$,

and similarly for $I^L, I^U$, and $F^L, F^U$ for the intervals $I$, and respectively $F$.

Hong-yu Zhang, Jian-qiang Wang, and Xiao-hong Chen [2] in 2014 defined the Interval Neutrosophic Score, Accuracy, and Certainty Functions as follows.

Let's consider *int([0, 1])* the set of all (open, semi-open/semi-closed, or closed) intervals included in or equal to *[0, 1]*, where the abbreviation and index *int* stand for *interval*, and *Zhang* stands for Hong-yu Zhang, Jian-qiang Wang, and Xiao-hong Chen.

### 9.1. Zhang Interval Neutrosophic Score Function

$$s_{int}^{Zhang} : \{int([0,1])\}^3 \to int([0,1])$$

$$s_{int}^{Zhang}(A) = [T^L + 1 - I^U + 1 - F^U, T^U + 1 - I^L + 1 - F^L]$$

### 9.2. Zhang Interval Neutrosophic Accuracy Function

$$a_{int}^{Zhang} : \{int([0,1])\}^3 \to int([0,1])$$

$$a_{int}^{Zhang}(A) = [min\{T^L - F^L, T^U - F^U\}, max\{T^L - F^L, T^U - F^U\}]$$

### 9.3. Zhang Interval Neutrosophic Certainty Function

$$c_{int}^{Zhang} : \{int([0,1])\}^3 \to int([0,1])$$

$$c_{int}^{Zhang}(A) = [T^L, T^U]$$

**10. New Interval Neutrosophic Score, Accuracy, and Certainty Functions**

Since comparing/ranking two intervals is equivalent to comparing/ranking two members (i.e. the intervals' midpoints), we simplify Zhang Interval Neutrosophic Score ($s_{int}^{Zhang}$), Accuracy

($a_{\text{int}}^{Zhang}$), Certainty ($c_{\text{int}}^{Zhang}$) functions, as follows:

$$s_{\text{int}}^{FS} : \{\text{int}([0,1])\}^3 \to [0,1]$$
$$a_{\text{int}}^{FS} : \{\text{int}([0,1])\}^3 \to [-1,1]$$
$$c_{\text{int}}^{FS} : \{\text{int}([0,1])\}^3 \to [0,1]$$

where the upper index *FS* stands for our name's initials, in order to distinguish these new functions from the previous ones:

### 10.1. New Interval Neutrosophic Score Function

$s_{\text{int}}^{FS}\left(([T^L, T^U], [I^L, I^U], [F^L, F^U])\right) = \frac{T^L + T^U + (1-I^L) + (1-I^U) + (1-F^L) + (1-F^U)}{6} = \frac{4 + T^L + T^U - I^L - I^U - F^L - F^U}{6}$,

which means the average of six positive-nesses;

### 10.2. New Interval Neutrosophic Accuracy Function

$a_{\text{int}}^{FS}\left(([T^L, T^U], [I^L, I^U], [F^L, F^U])\right) = \frac{T^L + T^U - F^L - F^U}{2}$, which means the average of differences between positiveness and negativeness;

### 10.3. New Interval Neutrosophic Certainty Function

$c_{\text{int}}^{FS}\left(([T^L, T^U], [I^L, I^U], [F^L, F^U])\right) = \frac{T^L + T^U}{2}$, which means the average of two positive-nesses.

### 10.4. Theorem

Let $\mathcal{M}_{int} = \{(T, I, F), \text{where } T, I, F \subseteq [0, 1], \text{ and } T, I, F \text{ are intervals}\}$, be the set of interval neutrosophic triplets.

The New Interval Neutrosophic Score, Accuracy, and Certainty Functions determine a total order relationship on the set $\mathcal{M}_{int}$ of Interval Neutrosophic Triplets.

*Proof*

Let's assume we have two interval neutrosophic triplets:
$P_1 = ([T_1^L, T_1^U], [I_1^L, I_1^U], [F_1^L, F_1^U])$,
and $P_1 = ([T_2^L, T_2^U], [I_2^L, I_2^U], [F_2^L, F_2^U])$, both from $M_{int}$.
We have to prove that: either $P_1 > P_2$, or $P_1 < P_2$, or $P_1 = P_2$.
Apply the new interval neutrosophic score function ($s_{\text{int}}^{FS}$) to both of them:

$$s_{\text{int}}^{FS}(P_1) = \frac{4 + T_1^L + T_1^U - I_1^L - I_1^U - F_1^L - F_1^U}{6}$$

$$s_{int}^{FS}(P_2) = \frac{4 + T_2^L + T_2^U - I_2^L - I_2^U - F_2^L - F_2^U}{6}$$

If $s_{int}^{FS}(P_1) > s_{int}^{FS}(P_2)$, then $P_1 > P_2$.

If $s_{int}^{FS}(P_1) < s_{int}^{FS}(P_2)$, then $P_1 < P_2$.

If $s_{int}^{FS}(P_1) = s_{int}^{FS}(P_2)$, then we get from equating the above two equalities that:

$T_1^L + T_1^U - I_1^L - I_1^U - F_1^L - F_1^U = T_2^L + T_2^U - I_2^L - I_2^U - F_2^L - F_2^U$

In this problematic case, we apply the new interval neutrosophic accuracy function ($a_{int}^{FS}$) to both $P_1$ and $P_2$, and we get:

$$a_{int}^{FS}(P_1) = \frac{T_1^L + T_1^U - F_1^L - F_1^U}{2}$$

$$a_{int}^{FS}(P_2) = \frac{T_2^L + T_2^U - F_2^L - F_2^U}{2}$$

If $a_{int}^{FS}(P_1) > a_{int}^{FS}(P_2)$, then $P_1 > P_2$.

If $a_{int}^{FS}(P_1) < a_{int}^{FS}(P_2)$, then $P_1 < P_2$.

If $a_{int}^{FS}(P_1) = a_{int}^{FS}(P_2)$, then we get from equating the two above equalities that:

$T_1^L + T_1^U - F_1^L - F_1^U = T_2^L + T_2^U - F_2^L - F_2^U$

Again, a problematic case, so we apply the new interval neutrosophic certainty function ($c_{int}^{FS}$) to both $P_1$ and $P_2$, and we get:

$$c_{int}^{FS}(P_1) = T_1^L + T_1^U$$

$$c_{int}^{FS}(P_2) = T_2^L + T_2^U$$

If $c_{int}^{FS}(P_1) > c_{int}^{FS}(P_2)$, then $P_1 > P_2$.

If $c_{int}^{FS}(P_1) < c_{int}^{FS}(P_2)$, then $P_1 < P_2$.

If $c_{int}^{FS}(P_1) = c_{int}^{FS}(P_2)$, then we get:

$T_1^L + T_1^U = T_2^L + T_2^U$

We prove that in the last case we get:

$P_1 =_N P_2$ (or $P_1$ is neutrosophically equal to $P_2$).

We get the following linear algebraic system of 3 equations and 12 variables:

$$\begin{cases} T_1^L+T_1^U - I_1^L-I_1^U - F_1^L-F_1^U = T_2^L+T_2^U - I_2^L-I_2^U - F_2^L-F_2^U \\ T_1^L+T_1^U - F_1^L-F_1^U = T_2^L+T_2^U - F_2^L-F_2^U \\ T_1^L+T_1^U = T_2^L+T_2^U \end{cases}$$

Second equation minus the third equation gives us:

$-F_1^L-F_1^U = -F_2^L-F_2^U$, or $F_1^L + F_1^U = F_2^L + F_2^U$.

First equation minus the second equation gives us:

$-I_1^L-I_1^U = -I_2^L-I_2^U$, or $I_1^L + I_1^U = I_2^L + I_2^U$.

The previous system is now equivalent to the below system:

$$\begin{cases} T_1^L+T_1^U = T_2^L+T_2^U \\ I_1^L+I_1^U = I_2^L+I_2^U \\ F_1^L+F_1^U = F_2^L+F_2^U \end{cases} \Leftrightarrow \begin{cases} \dfrac{T_1^L+T_1^U}{2} = \dfrac{T_2^L+T_2^U}{2} \\ \dfrac{I_1^L+I_1^U}{2} = \dfrac{I_2^L+I_2^U}{2} \\ \dfrac{F_1^L+F_1^U}{2} = \dfrac{F_2^L+F_2^U}{2} \end{cases}$$

which means that:

*i)* the intervals $[T_1^L, T_1^U]$ and $[T_2^L, T_2^U]$ have the same midpoint, therefore they are neutrosophically equal.

*ii)* the intervals $[I_1^L, I_1^U]$ and $[I_2^L, I_2^U]$ have also the same midpoint, so they are neutrosophically equal.

*iii)* similarly, the intervals $[F_1^L, F_1^U]$ and $[F_2^L, F_2^U]$ have the same midpoint, and again they are neutrosophically equal.

Whence, the interval neutrosophic triplets $P_1$ and $P_2$ are neutrosophically equal, i.e. $P_1 =_N P_2$.

**10.5. Theorem**

Let's consider the ranking of intervals defined by Xu and Da, which is equivalent to the ranking of intervals' midpoints. Then, the algorithm by Hong-yu Zhang et al. for ranking the interval neutrosophic triplets in equivalent to our algorithm.

*Proof*

Let's consider two interval neutrosophic triplets, $P_1$ and $P_2 \in \mathcal{M}_{int}$,

$P_1 = ([T_1^L, T_1^U], [I_1^L, I_1^U], [F_1^L, F_1^U])$,

and $P_2 = ([T_2^L, T_2^U], [I_2^L, I_2^U], [F_2^L, F_2^U])$.

Let's rank them using both methods and prove we get the same results. We denote by $s_{\text{int}}^{Zhang}$, $a_{\text{int}}^{Zhang}$, $c_{\text{int}}^{Zhang}$, and $s_{int}^{FS}$, $a_{int}^{FS}$, $c_{int}^{FS}$ the Interval Neutrosophic Score, Accuracy, and Certainty Functions, by Hong-yu Zhang et al. and respectively by us.

*Interval Neutrosophic Score Function*

$$s_{\text{int}}^{Zhang}(P_1) = [T_1^L + 1 - I_1^U + 1 - F_1^U, T_1^U + 1 - I_1^L + 1 - F_1^L]$$

$$s_{\text{int}}^{Zhang}(P_2) = [T_2^L + 1 - I_2^U + 1 - F_2^U, T_2^U + 1 - I_2^L + 1 - F_2^L]$$

a) If $s_{\text{int}}^{Zhang}(P_1) > s_{\text{int}}^{Zhang}(P_2)$, then

the midpoint of the interval $s_{\text{int}}^{Zhang}(P_1) >$ midpoint of the interval $s_{\text{int}}^{Zhang}(P_2)$,

or $\frac{T_1^L + 1 - I_1^U + 1 - F_1^U + T_1^U + 1 - I_1^L + 1 - F_1^L}{2} > \frac{T_2^L + 1 - I_2^U + 1 - F_2^U + T_2^U + 1 - I_2^L + 1 - F_2^L}{2}$,

or $T_1^L + T_1^U - I_1^L - I_1^U - F_1^L - F_1^U > T_2^L + T_2^U - I_2^L - I_2^U - F_2^L - F_2^U$,

or $\frac{4 + T_1^L + T_1^U - I_1^L - I_1^U - F_1^L - F_1^U}{6} > \frac{4 + T_2^L + T_2^U - I_2^L - I_2^U - F_2^L - F_2^U}{6}$,

or $s_{int}^{FS}(P_1) > s_{int}^{FS}(P_2)$.

b) If $s_{\text{int}}^{Zhang}(P_1) < s_{\text{int}}^{Zhang}(P_2)$, the proof is similar, we only replace the inequality symbol $>$ by $<$ into the above proof.

c) If $s_{\text{int}}^{Zhang}(P_1) = s_{\text{int}}^{Zhang}(P_2)$, the proof again is similar with the above, we only replace $>$ by $=$ into the above proof.

*Interval Neutrosophic Accuracy Function*

$$a_{\text{int}}^{Zhang}(P_1) = [\min\{T_1^L - F_1^L, T_1^U - F_1^U\}, \max\{T_1^L - F_1^L, T_1^U - F_1^U\}];$$

$$a_{\text{int}}^{Zhang}(P_2) = [\min\{T_2^L - F_2^L, T_2^U - F_2^U\}, \max\{T_2^L - F_2^L, T_2^U - F_2^U\}].$$

a) If $a_{\text{int}}^{Zhang}(P_1) > a_{\text{int}}^{Zhang}(P_2)$, then

the midpoint of the interval $a_{int}^{Hong}(P_1) >$ the midpoint of the interval $a_{int}^{Hong}(P_2)$,

or $\frac{T_1^L - F_1^L + T_1^U - F_1^U}{2} > \frac{T_2^L - F_2^L + T_2^U - F_2^U}{2}$,

or $a_{int}^{FS}(P_1) > a_{int}^{FS}(P_2)$.

b) Similarly, if $a_{\text{int}}^{Zhang}(P_1) < a_{\text{int}}^{Zhang}(P_2)$, just replacing $>$ by $<$ into the above proof.

c) Again, similarly if $a_{\text{int}}^{Zhang}(P_1) = a_{\text{int}}^{Zhang}(P_2)$, only replacing $>$ by $=$ into the above proof.

*Interval Neutrosophic Certainty Function*

$$c_{\text{int}}^{Zhang}(P_1) = [T_1^L, T_1^U]$$

$$c_{\text{int}}^{Zhang}(P_2) = [T_2^L, T_2^U]$$

a) If $c_{\text{int}}^{Zhang}(P_1) > c_{\text{int}}^{Zhang}(P_2)$, then

b) **the midpoint of the interval** $c_{\text{int}}^{Zhang}(P_1) >$ **the midpoint of the interval** $c_{\text{int}}^{Zhang}(P_2)$, or $\frac{T_1^L + T_1^U}{2} > \frac{T_2^L + T_2^U}{2}$ or $c_{int}^{FS}(P_1) > c_{int}^{FS}(P_2)$.

b) Similarly, if $c_{\text{int}}^{Zhang}(P_1) < c_{\text{int}}^{Zhang}(P_2)$, just replacing > by < into the above proof.

c) Again, if $c_{\text{int}}^{Zhang}(P_1) = c_{\text{int}}^{Zhang}(P_2)$, only replace > by = into the above proof.

Therefore, we proved that, for any interval neutrosophic triplet $P$,

$s_{\text{int}}^{Zhang}(P) \sim s_{int}^{FS}(P)$, where ~ means equivalent;

$a_{\text{int}}^{Zhang}(P) \sim a_{int}^{FS}(P)$,

and $c_{\text{int}}^{Zhang}(P) \sim c_{int}^{FS}(P)$.

## 11. Subset Neutrosophic Score, Accuracy, and Certainty Functions

Let $M_{\text{subset}} = \{(T_{\text{subset}}, I_{\text{subset}}, F_{\text{subset}}), \text{where the subsets } T_{\text{subset}}, I_{\text{subset}}, F_{\text{subset}} \subseteq [0,1]\}$.

We approximate each subset by the smallest closed interval that includes it.

Let's denote:

$T^L = inf(T_{\text{subset}})$ and $T^U = sup(T_{\text{subset}})$; therefore $T_{\text{subset}} \subseteq [T^L, T^U]$;

$I^L = inf(I_{\text{subset}})$ and $I^U = sup(I_{\text{subset}})$; therefore $I_{\text{subset}} \subseteq [I^L, I^U]$;

$F^L = inf(F_{\text{subset}})$ and $F^U = sup(F_{\text{subset}})$; therefore $F_{\text{subset}} \subseteq [F^L, F^U]$.

Then:

$$M_{\text{subset}} \approx \left\{ \begin{array}{l} ([T^L, T^U], [I^L, I^U], [F^L, F^U]), \text{where } T^L, T^U, I^L, I^U, F^L, F^U \in [0,1], \\ \text{and } T^L \leq T^U, I^L \leq I^U, F^L \leq F^U \end{array} \right\}$$

### 11.1. Definition of Subset Neutrosophic Score, Accuracy, and Certainty Functions

Then, the formulas for Subset Neutrosophic Score, Accuracy, and Certainty Functions will coincide with those for Interval Neutrosophic Score Accuracy, and Certainty Functions by Hong-yu Zhang, and respectively by us.

### 11.2. Theorem

Let N be the Interval Neutrosophic Triplet

$N = ([T^L, T^U], [I^L, I^U], [F^L, F^U])$,

where $T^L \leq T^U$, $I^L \leq I^U$, $F^L \leq F^U$,

and all $[T^L, T^U], [I^L, I^U], [F^L, F^U] \subseteq [0, 1]$.

If each interval collapses to a single point, i.e.

$T^L = T^U = T$, then $[T^L, T^U] = [T, T] \equiv T \in [0, 1]$,

$I^L = I^U = I$, then $[I^L, I^U] = [I, I] \equiv I \in [0, 1]$,

$F^L = F^U = F$, then $[F^L, F^U] = [F, F] \equiv T \in [0, 1]$,

then $s_{int}^{FS}(N) = s(N)$,

$a_{int}^{FS}(N) = a(N)$,

and $c_{int}^{FS}(N) = c(N)$.

*Proof*

$$s_{int}^{FS}(N) = \frac{4 + T^L + T^U - I^L - I^U - F^L - F^U}{6} = \frac{4 + T + T - I - I - F - F}{6}$$

$$= \frac{4 + 2T - 2I - 2F}{6} = \frac{2 + T - I - F}{3} = s(N).$$

$$a_{int}^{FS}(N) = \frac{T^L + T^U - F^L - F^U}{2} = \frac{T + T - F - F}{2} = \frac{2(T - F)}{2} = a(N).$$

$$c_{int}^{FS}(N) = \frac{T^L + T^U}{2} = \frac{T + T}{2} = \frac{2T}{2} = c(N).$$

### 12. Conclusion

The most used and easy for ranking the Neutrosophic Triplets (*T, I, F*) are the following functions, that provide a total order:

**Single-Valued Neutrosophic Score, Accuracy, and Certainty Functions:**

$$s(T, I, F) = \frac{2 + T - I - F}{3}$$

$$a(T, I, F) = T - F$$

$$c(T, I, F) = T$$

**Interval-Valued Neutrosophic Score, Accuracy, and Certainty Functions:**

$$s_{\text{int}}^{FS}\left(([T^L, T^U], [I^L, I^U], [F^L, F^U])\right) = \frac{4 + T^L + T^U - I^L - I^U - F^L - F^U}{6}$$

$$a_{\text{int}}^{FS}\left(([T^L, T^U], [I^L, I^U], [F^L, F^U])\right) = \frac{T^L + T^U - F^L - F^U}{2}$$

$$c_{\text{int}}^{FS}\left(([T^L, T^U], [I^L, I^U], [F^L, F^U])\right) = \frac{T^L + T^U}{2}$$

All these functions are very much used in decision-making applications.

**References**

[1] Florentin Smarandache, *The Average Positive Qualitative Neutrosophic Function and The Average Negative Qualitative Neutrosophic Function*, section I.2, pp. 20-24, in his book *Neutrosophic Perspectives: Triplets, Duplets, Multisets, Hybrid Operators, Modal Logic, Hedge Algebras. And Applications (Second extended and improved edition),* Pons Publishing House, Brussels, 348 p., 2017, http://fs.unm.edu/NeutrosophicPerspectives-ed2.pdf

[2] Z. S. Xu and Q. L. Da, *The uncertain OWA operator*, International Journal of Intelligent Systems, 569-575, Vol. 17, 2002.

[3] Hong-yu Zhang, Jian-qiang Wang, and Xiao-hong Chen, *Interval Neutrosophic Sets and Their Application in Multicriteria Decision Making Problems*, Hindawi Publishing Corporation, The Scientific World Journal, Volume 2014, 15 p., *http://dx.doi.org/10.1155/2014/64595*

[4] Jiqian Chen and Jun Ye, *Some Single-Valued Neutrosophic Dombi Weighted Aggregation Operators for Multiple Attribute Decision-Making*, Licensee MDPI, Basel, Switzerland, Symmetry, 9, 82, 2017; DOI:10.3390/sym9060082

[5] Jun Ye, Communications by emails to the author, 10-11 October 2020.